\documentclass[]{article}%elsart

\usepackage{amsmath}
\usepackage{amssymb}
\usepackage[ansinew]{inputenc}

\def\Nset{\mathbb{N}}

\def\Rset{\mathbb{R}}

\newtheorem{theorem}{Theorem}
\newtheorem{lemma}{Lemma}
\newtheorem{example}{Example}
\newtheorem{cor}{Corollary}

\newtheorem{remark}{Remark}
\newtheorem{prop}{Proposition}
\newcommand{\proof}{\bfseries Proof: $\quad$\mdseries}
\newcommand{\finishproof}{\begin{flushright} $\Box$ \end{flushright}}

\date{\today}

\begin{document}

%\begin{frontmatter}
\title{$L^{2}$-spectral gaps for time discrete reversible Markov chains}
\author{Achim Wübker\thanks{e-mail: awuebker@mathematik.uni-osnabrueck.de}\\Institute of Mathematics, University of Osnabrück\\Albrechtstraße 28 a, 49076 Osnabrück}

\maketitle

\begin{abstract}
In this paper we study the spectral properties of Markov-operator on $L^{2}$-spaces.
Lawler and Sokal (Trans. Amer. Math. Soc., 1988, 309, pp. 557-580) used isoperimetric constants for discrete and continuous time Markov chains to obtain a spectral gap at 1. For time discrete Markov chains
this does not exclude periodic behavior. We define a new constant measuring the distance from periodicity and give necessary and sufficient conditions for the existence of a global spectral gap in terms of this constant. 
\end{abstract}

\section{Introduction}
Several techniques had been developed for studying the stochastic behavior of Markov processes. For example, one may analyze the spectral properties of the associated Markov operator. A well-known result is the 
Frobenius-Perron theorem, which states that for an irreducible, aperiodic and finite state space Markov chain $\xi_1,\xi_2,\ldots$ we get that $1$ is an eigenvalue (according to a left eigenvector $\pi$, $\pi$ is the uniquely determined invariant measure) of $P$ ($P$ the transition matrix of the chain) and the absolute values of all other eigenvalues $\lambda$ are smaller than one. Moreover, the speed of convergence of the transition probabilities to  $\pi$ is geometric and can be estimated in terms of $|\lambda_2|$, where $\lambda_2$ is the eigenvalue with second largest modulus. The quantity $r:=1-|\lambda_2|$ is often called the spectral gap of the chain. Since for Markov chains with large state space it is often impossible to compute $\lambda_2$, one is interested in estimating the gap. In the 
last 20 years several papers concerning spectral gaps for Markov chains has been published (e.g. \cite{chen}, \cite{chenbuch}, \cite{chenwang}, \cite{diaconis3}, \cite{diaconis4}, \cite{diaconis1}, \cite{diaconis2}, \cite{fill}, \cite{jerrum}, \cite{lawler}). There are mainly two reasons for this: On one hand, in the early 1990`s 
Markov Chain Monte Carlo (MCMC) methods were developed and statements concerning the speed of convergence needed information about the size of the spectral gap. On the other hand, at that time the mathematicians
working in this area started to use the technique of Dirichlet Forms which proved to be very successful
to tackle this problem. 
%Most papers in the literature treated finite or countable state space MC
%(e.g. \cite{diaconis1}, \cite{diaconis2}, \cite{diaconis3}, \cite{diaconis4}, \cite{fill}). Only a 
%few papers are concerned with general state space MC. In this paper we get estimates for use the approach 
%Lawler and Sokal \cite{lawler} using isoperimetric constants. \\    

First, Cheeger \cite{cheeger} introduced an isoperimetric constant $k$ of a
compact Riemannian manifold $M$ in order to give a lower bound for the smallest strictly positive
eigenvalue of the Laplacian on $M$. Since a discrete version of the Laplacian can be associated with time discrete Markov chains, several authors investigated spectral properties of Markov operators using
Cheeger's idea (e.g. \cite{chen}, \cite{chenbuch}, \cite{chenwang}, \cite{dodziuk}, \cite{fill},  \cite{lawler}). A fundamental paper comprising the results of research in this area until 1988 is provided by Lawler and Sokal \cite{lawler}.\\
Lawler and Sokal \cite{lawler} established a simple necessary and sufficient condition for reversible, continuous time Markov chains that ensures the existence of a spectral gap in terms of isoperimetric constant $k$. Their estimates concerning the size of the gap were partially improved by Chen \cite{chen} in 2000. The condition mentioned before is still necessary for discrete time Markov chain, but it is not sufficient. As far as we know, in this setting no simple condition equivalent to the spectral gap property has appeared in the literature. 
We introduce a constant $K$ which measures the chain's distance from
periodicity. This constant is used to obtain a condition that ensures the spectral gap property for 
reversible discrete time Markov chain on a general state space. Moreover, we establish a lower bound for the size of the gap in terms of $k$ and $K$. For this, we use a geometric approach which can be generalized to the non-reversible case.

\section{A spectral gap theorem for reversible Markov chains}
In this section we introduce the basic notations and remind some well known facts that will be needed in the sequel. Throughout this paper we consider a positive recurrent Markov chain $\xi_1,\xi_2,\ldots$ with arbitrary state space $(\Omega,\mathcal{F})$, transition kernel $p(x,dy)$ and uniquely determined invariant probability measure $\pi$. 
The stochastic kernel $p(x,dy)$ induces a linear operator $P$ on the Hilbert space $L^{2}(\pi)$ by
\begin{equation}
Pf(x):=\int_{\Omega}f(y)p(x,dy).  
\end{equation}
Using Jensen's inequality, we see that $||P||_{2}\le 1$ and since $1$ is an eigenfunction of $P$, we have
$||P||_{2}=1$. 
Sometimes we will assume that the Markov chain is reversible, i.e.
\begin{equation}
\pi(dx)p(x,dy)=\pi(dy)p(y,dx).
\end{equation}
This implies that $P$ is self adjoint and the spectrum of the operator $P$ is real valued, i.e. \\$\sigma(P)\subset[-1,1]$.\\

The main question of interest is now whether an $r>0$ exists with \begin{equation}\label{spectral gap}
\lim_{n\rightarrow\infty}||P^{n}f||_{2}^{\frac{1}{n}}\le 1-r.
\end{equation}
for all $f\in L_{0,1}^{2}(\pi):=\{f\in L^{2}(\pi):\int_{\Omega}f(x)\pi(dx)=0,\int_{\Omega}f(x)^{2}\pi(dx)=1\}$. 
The supremum over all $r>0$ satisfying (\ref{spectral gap}) is called the spectral gap of $P$.\\
In order to estimate $r$, let us define
\begin{equation}
k:=\inf_{A\in\mathcal{F}:0<\pi(A)<\frac{1}{2}}k(A),
\end{equation}
as introduced in \cite{lawler}, where 
\begin{equation}
k(A):=\frac{1}{\pi(A)\pi(A^{c})}\int_{A}p(x,A^{c})\pi(dx).
\end{equation}
The quantity $k(A)$ measures the probability flow out of the set $A$ to its complement.\\ 

Let us introduce a family $k_n$ of isoperimetric constants associated with the Markov chain $\xi_1,\xi_{n+1},\xi_{2n+1},\ldots$.
\begin{equation}\label{isoneu}
k_{n}:=\inf_{A\in\mathcal{F}}k_{n}(A),\quad k_{n}(A):=\frac{1}{\pi(A)\pi(A^{c})}\int_{A}p^{n}(x,A^{c})\pi(dx), \,\,n\in\Nset.
\end{equation}
From the definition it follows that $k=k_{1}$.\\

Now let us assume that $P$ is self adjoint. 
We say that the operator $P$ has a spectral gap of size $r_{1}$ at 1 if and only if  
\begin{equation}\label{atone}
\sigma(P)\subset[-1,1-r_{1}],\,\,\,\,r_1>0.
\end{equation}
Accordingly, we say that $P$ has a spectral gap at $-1$ if and only if a $r_{-1}>0$  exists such that
\begin{equation}\label{atminusone}
\sigma(P)\subset[-1+r_{-1},1]
\end{equation}
The size of the gap is given by the supremum over all $r_{1}$ ($r_{-1}$) such that (\ref{atone}) ((\ref{atminusone}) respectively) is fulfilled.
Lawler and Sokal \cite{lawler} proved the following bounds for the spectral gap at 1 in terms of $k$:
\begin{theorem}[Lawler, Sokal]\label{LS}
For a reversible, positive recurrent Markov chain with uniquely determined stationary distribution $\pi$
we know that
\begin{equation}\label{ELS}
k\ge r_{1}\ge\frac{\kappa}{8}k^{2},
\end{equation}
\end{theorem}
where $\kappa$ is a constant larger than one (see \cite{lawler}). 

For reversible finite state space Markov chain a similar result is provided by Jerrum and Sinclair \cite{jerrum} in terms of conductance (for the definition of conductance see e.g. \cite{bremaud}), which is closely related to the definition of the isoperimetric constant.

In order to have a spectral gap $r$ in the sense of (\ref{spectral gap}) we should also know 
the spectral gap $r_{-1}$ at $-1$, since $r=\min(r_{-1},r_{1})$.
The value of $r_{-1}$ measures the periodic behavior of the chain in a certain sense. For this
reason let us consider the following \vspace{0.5cm}\\

\begin{example}
Let $\Omega:=\{1,2,3,4\}$. Assume that\\
$p(i,i+1(\,mod\,4))=p(i,i-1(\,mod\,4))=\frac{1}{2}$, $i\in\Omega$. Then the invariant measure
is given by 
$\pi=(\frac{1}{4},\frac{1}{4},\frac{1}{4},\frac{1}{4})$. Let us consider the set
$A=\{1,3\}$ and define $f:=1_{A}-1_{A^{c}}$.

One can show that $k>0$ and
\[Pf=-f,\]
so $-1$ is an eigenvalue of $P$. Hence this chain does not have a spectral gap in 
the sense of (\ref{spectral gap}).\vspace{0.2cm}\\
%As we will see, the constants $k(A)$ can be used to measure the periodic behavior of the chain, since the probability flow in the periodic case would be large in a sense which will be more specified now. 

\end{example}

%Especially $\epsilon_{1}=0$ means that the chain is more of less periodic of order two. Higher order periodicities are excluded by the reversibility of the chain.
%It is easy to see that $P$ is selfadjoined on $L^{2}(\pi)$ and $-1\in\sigma(P)$.
%This is equivalent to $\epsilon_{-1}=0$. 
%For this reason let us return to the example above. 
The following observation is crucial. In the example above we have
\[k(A)=2.\]
Let us define a new constant, namely
\[K:=\sup_{A\in\mathcal{F}}k(A).\] 
In the way $k$ measures the distance from ergodicity, $2-K$ measures the distance from 
periodicity. Since periodicity of order larger than 2 is excluded by the reversibility assumption, one may
hope that $K<2$ is actually a sufficient condition. \vspace{0.5cm} \\

We now state the main result of the paper, namely a condition that ensures the existence of a spectral gap for general state space Markov chains in terms
of $k$, $k_2$ and $K$.\\
 
\begin{theorem}\label{mainresult}
Let $\xi_1,\xi_2,\ldots$ be a reversible, time homogeneous Markov chain with stationary distribution $\pi$ and transition kernel $p(\cdot,\cdot)$.  
Then the following three conditions are equivalent:
\begin{enumerate}
\item
$P$ has a $L^{2}$-spectral gap.
\item
\begin{equation}\label{general}
0<k\le K<2.
\end{equation}
\item
\begin{equation}
k_2>0.
\end{equation}
\end{enumerate}
For the spectral gap we obtain the following estimate:
\begin {equation}
\sigma(P)\subset\left[-\sqrt{1-\frac{\kappa}{8}k_{2}^{2}},\min(\sqrt{1-\frac{\kappa}{8}k_{2}^{2}},1-\frac{\kappa}{8}k^2)\right],
\end{equation}
where $k_2$ can be estimated from below by

\begin{eqnarray}
k_2&\ge&\sup_{\delta,\epsilon_{1},\epsilon_{2},\epsilon\in\Rset_{+}}\min \left[\frac{k^{2}}{16}\delta,\frac{k}{4}(
\epsilon_{1}\epsilon_{2}(1-\delta)-\delta),\right.\nonumber\\
&&\,\,\,\,\,\,\,\,\left.\left(k\left(\frac{(2-\epsilon)(1-\epsilon_1)(1-\epsilon_2)(1-\delta)}{(1-\epsilon)K}-\frac{1}{1-\epsilon}\right)-\frac{\epsilon}{1-\epsilon}\right)\epsilon\right].
\end{eqnarray}
\end{theorem}
It is known that the first inequality is a necessary condition for spectral gap property (\cite{jerrum}, \cite{lawler}).
For continuous time Markov chains, the condition $\tilde{k}>0$ ($\tilde{k}$ defined by the transition rates $q(\cdot,\cdot)$ instead of the transition probabilities $p(\cdot,\cdot)$) is indeed sufficient, since periodicity is a phenomenon occurring only for discrete time Markov chains. As mentioned above, the third inequality in (\ref{general}) excludes periodicity of the chain.

\section{Proof of Theorem 2}
In order to prove the theorem above, one has to establish several lemmas.
%, some of them are probably 
%well known. Since we didn't found references in the literature, we present also the proofs. 
In the following we use the setting above. If reversibility is required, then it will be mentioned explicitly.

\begin{lemma}\label{first}
For all $A\in\mathcal{F}$ holds:
\[k(A)=k(A^{c}).\]
\end{lemma}
\proof
We see that for the proof we only need the stationarity of the starting distribution $\pi$ of the Markov chain.
\begin{eqnarray}
\int_{A}p(x,A^{c})\pi(dx)&=&\int_{\Omega}p(x,A^{c})\pi(dx)-\int_{A^{c}}p(x,A^{c})\pi(dx)\nonumber\\
&=&\pi(A^{c})-(\pi(A^{c})-\int_{A^{c}}p(x,A)\pi(dx))\nonumber\\
&=&\int_{A^{c}}p(x,A)\pi(dx).\nonumber
\end{eqnarray}
\finishproof

The next lemma provides a precise upper bound for $k(A)$:

\begin{lemma}\label{secondone}
For all $A\in\mathcal{F}$ we have
\begin{equation} 
0\le k(A)\le 2.
\end{equation}
\end{lemma}
\proof
The first inequality is trivial. In order to show the second one, let us consider $K:=\sup_{A\in\mathcal{F}}k(A)$.
Then for every $\epsilon>0$ there exists $A\in\mathcal{F}$ such that $k(A)>K-\epsilon$. For such a set $A$ it follows

\begin{equation}\label{notd}
\int_{A}p(x,A^{c})\pi(dx)>(K-\epsilon)\pi(A)\pi(A^{c}).
\end{equation}
Since $p(x,A^{c})$ only takes values between $0$ and $1$, an $\alpha$ with $0\le\alpha\le 1$ exists such that
\begin{equation}\label{imfolbe}
\int_{A}p(x,A^{c})\pi(dx)=\alpha \pi(A).
\end{equation}
This together with (\ref{notd}) yields
\begin{equation}\label{lutz5}
\alpha>(K-\epsilon)\pi(A^{c}).
\end{equation}
By Lemma \ref{first} we also have
\[\int_{A^{c}}p(x,A)\pi(dx)> (K-\epsilon)\pi(A)\pi(A^{c}).\]
So there exists $\beta\in[0,1]$ such that
\begin{equation}\label{imfolbe2}
\int_{A^{c}}p(x,A)\pi(dx)=\beta\pi(A^{c}).
\end{equation}
As above this yields
\begin{equation}\label{lutz6}
\beta > (K-\epsilon)\pi(A).
\end{equation}
Adding (\ref{lutz5}) and (\ref{lutz6}) we get
\begin{equation}\label{lutz10}
2\ge \alpha +\beta > (K-\epsilon) \pi(A) +(K-\epsilon)\pi(A^{c})=K-\epsilon.
\end{equation}
Since this is true for all $\epsilon>0$, we have $K\le 2$, which proves the lemma.
\finishproof

We saw in Example 1 that $K=2$ is possible and the bound was attained by some set $A$ with $\pi(A)=\frac{1}{2}$.
The next lemma shows that $k(A)>2-\epsilon$ is only possible for sets $A\in\mathcal{F}$ with $\pi(A)$ close to $\frac{1}{2}$.

\begin{lemma}\label{einhalb}
Let $K:=\sup_{A\in\mathcal{F}}k(A)=2$. Then for any sequence of sets $A_n\in\mathcal{F}$ with the property that $k(A_n)\rightarrow 2$, we have that $\pi(A_n)\rightarrow\frac{1}{2}$.
\end{lemma}
\proof  
According to (\ref{imfolbe}) and (\ref{imfolbe2}), we associate with $A_n$ and $A_{n}^{c}$ the constants $\alpha_{n}$ and $\beta_n$. From $k(A_n)\ge 2-\epsilon_n$, $0<\epsilon_n\rightarrow 0$ and (\ref{lutz10}) with $K=2$, $\alpha=\alpha_n$, $\beta=\beta_n$, $\epsilon=\epsilon_n$  it follows that $\alpha_n\rightarrow 1$ and $\beta_n\rightarrow 1$. Using (\ref{lutz5}) and (\ref{lutz6}) with $K=2$
and keeping in mind the first inequality in (\ref{lutz10}) we can conclude that $\pi(A_n)\rightarrow\frac{1}{2}$.
\finishproof

We have the following simple relation between $k$ and $k_n$: %as defined in \ref{isoneu} 
\begin{lemma}\label{Potenzen}
\[k=0\Longleftrightarrow\, k_n=0 \quad\forall n\in\Nset.\]
%then
%\[k_{n}:=\inf_{A\in\mathcal{F}}k_{n}(A)=0\,\,\forall n\in\Nset,\]
%wobei
%$k_n(A):=\frac{1}{\pi(A)\pi(A^{c})}\int_{A}p^{n}(x,A^{c})\pi(dx)$.
\end{lemma}
\proof
\begin{eqnarray}
\pi(A)\pi(A^{c})k_n(A)&=&\int_{A}p^{n}(x,A^{c})\pi(dx)=\int_{A}\int_{A}p(x,dy)p^{n-1}(y,A^{c})\pi(dx)\nonumber\\
&+&\int_{A}\pi(dx)\int_{A^{c}}p(x,dy)p^{n-1}(y,A^{c})\nonumber\\
&\le&(k_{n-1}(A)+k(A))\pi(A)\pi(A^{c})\le\ldots\le n\,k(A)\pi(A)\pi(A^{c})\nonumber.
\end{eqnarray}
\finishproof

%\begin{lemma}\label{dasbraucheichnoch}
%Assume that there exists $n_{0}\in\Nset,\epsilon>0$, such that $k_n\ge\epsilon\,\forall n\ge n_{0}$. %Then we have that
%\[k_i\ge\frac{\epsilon}{n_0}\,\forall i\in\Nset.\]
%\end{lemma}
%\proof
%For $k_1$ folgt die Absch"atzung direkt aus dem Beweis des vorherigen Lemmas. F"ur die "ubrigen 
%$k_i$, $i\in\{2,\ldots,n_{0}-1\}$ kann man das so sehen:
%\[\int_{A}p^{i}(x,A^{c})\ge \frac{1}{2}\int_{A}p^{2i}(x,A^{c})\ge\ldots \ge\frac{1}{2^{[\log_{2} %n_0]}}\int_{A}p^{i2^{[\log_{2} n_0]}}(x,A^{c}).\]
%Da 
%\[i2^{[\log_{2} n_0]}\ge n_{0}\,\,\forall i\in\{2,\ldots,n_{0}-1\}\] 
%und 
%\[2^{[\log_{2} n_0]}\le n_{0},\] 
%folgt die Behauptung.
%\finishproof

\begin{prop}\label{absch}
Let $P$:  $L_{0}^{2}(\pi)\rightarrow L_{0}^{2}(\pi)$ be the self adjoint operator associated to the reversible Markov chain $\xi_1,\xi_2,\ldots$. Then we obtain
\begin{equation}
\sigma(P)\subset\left[-\sqrt{1-\frac{\kappa}{8}k_{2}^{2}},\sqrt{1-\frac{\kappa}{8}k_{2}^{2}}\right].
\end{equation}
\end{prop}
\proof
Since $P$ is self adjoint, $P^{2}$ is self adjoint and positive. Applying Theorem 1 to $P^2$ yields
\[\sigma(P^2)\subset[0,1-\frac{\kappa}{8}k_{2}^{2}].\]
Using the spectral mapping theorem ( \cite{werner}), 
%\[\sigma((P-P_1)^{2})=(\sigma(P-P_1))^2,\]
we obtain 
\[\sigma(P)\subset\left[-\sqrt{1-\frac{\kappa}{8}k_{2}^{2}},\sqrt{1-\frac{\kappa}{8}k_{2}^{2}}\right].\]
\finishproof

\begin{cor}\label{majaja}
The Markov chain $\xi_1,\xi_2,\ldots$ has a spectral gap if and only if
\begin{equation}\label{preplenkel}
k_2>0.
\end{equation}
\end{cor}
\proof
The sufficiency of the condition follows immediately from Proposition \ref{absch}.
Now let us assume that $k_2=0$. Applying Lemma \ref{Potenzen} to the Markov chain
$\xi_1,\xi_3,\xi_5,\ldots$, we obtain that $k_{2n}=0$ for all $n\in\Nset$. Then applying
the first inequality in (\ref{ELS}) to the chain $\xi_1,\xi_{2n+1},\xi_{4n+1},\ldots$
yields 
\[\{1\}\subset\sigma(P^{2n})\mbox{ for all }n\in\Nset.\]
This finalizes the proof.
\finishproof
 
In order to prove Theorem 2, we have to establish the connection between $k,K$ on one
hand and $k_2$ on the other hand.

\begin{lemma}\label{leema}
Assume that $K=2$. Then $k_2=0$.
\end{lemma}
\proof
For each $A_n$ in $\mathcal{F}$ we obtain
\begin{eqnarray}
\pi(A_n)\pi(A_{n}^{c})k_{2}(A_n)&=&\int_{A_n}\int_{A_n}\pi(dx)p(x,dy)p(y,A_{n}^{c})\nonumber\\
&+&\int_{A_n}\int_{A_{n}^{c}}\pi(dx)p(x,dy)p(y,A_{n}^{c})\nonumber\\
&\le&\int_{A_n}\pi(dx)p(x,A_n)+\int_{A_{n}^{c}}\pi(dx)p(x,A_{n}^{c})\nonumber\\
&=&\pi(A_n)-\int_{A_{n}^{c}}\pi(dx)p(x,A_n)+\pi(A_{n}^{c})\nonumber\\
&&-\int_{A_n}\pi(dx)p(x,A_{n}^{c})\nonumber\\
&=& 1-2\int_{A_n}\pi(dx)p(x,A_{n}^{c}).\nonumber
\end{eqnarray}
Now choose a sequence $A_n\in\mathcal{F}$ with $k(A_n)\rightarrow 2$. From Lemma \ref{einhalb} we know that $\pi(A_n)\rightarrow\frac{1}{2}$. This yields
\[k_{2}(A_n)\le \frac{1}{\pi(A_n)\pi(A_{n}^{c})}-2k(A_n)\rightarrow 0,\,\,n\rightarrow\infty.\]
\finishproof

\begin{remark}
Since we did not use the reversibility assumption, Lemma \ref{leema} is true for arbitrary
Markov chains.
\end{remark}
Now we come to the most difficult part of the proof of Theorem 2:
%\begin{lemma}\label{tatbes}
%Es gelten wieder die zu Beginn des Abschnitts \ref{sectioniso1} gemachten Bedingungen. Sei dar"uber hinaus $(\Omega,\mathcal{F},\pi)$ ein Ma"sraum, sodass $\inf_{A\in\mathcal{F}:\pi(A)>0}\pi(A)=0$. Dann gilt:
%\[k\le 1.\]
%\beweis
%Nach Voraussetzung existiert eine Folge $A_n\in\mathcal{F}$ mit $\pi(A_n)\rightarrow0$. Dann gilt offensichtlich:
%\[k(A_n)=\frac{1}{\pi(A_{n}^{c})}\int_{A_n}p(x,A_{n}^{c})\frac{\pi(dx)}{\pi(A_n)}\le\frac{1}{\pi(A_{n}^{c})}.\]
%Limesbildung liefert die Behauptung.
%\beweisende
%\end{lemma}

%\begin{satz}\label{mathis}
%Sei $\xi_1, \xi_2, \ldots$ eine reversible MK mit station"arem Wma"s $\pi$, "Ubergangskern %$p(\cdot,\cdot)$ und $P$ der assoziierte Markov-Operator auf $L_2(\pi)$. Dann sind folgende Aussagen "aquivalent:
%\begin{itemize}
%\item
%\[0<k=\inf_{A\in\mathcal{F}}k(A)\le \sup_{A\in\mathcal{F}}k(A)<2.\]
%\item
%\[\exists C>0,\,\, q<1, \mbox{ sodass } ||(P-P_1)^{n}||_{2}\le C\,q^{n}\,\,\,\forall n\in\Nset.\]
%\end{itemize}
%\end{satz}

%Bevor dieser Satz bewiesen wird, ben"otigen wir noch ein technisches lemma:
\begin{lemma}
Let $\xi_1, \xi_2, \ldots$ be a reversible, time homogeneous Markov chain with arbitrary state space $(\Omega,\mathcal{F},\pi)$. We get the following estimate for $k_2$:
\begin{eqnarray*}\label{jakob7}
k_2&\ge&\sup_{\delta,\epsilon_{1},\epsilon_{2},\epsilon\in\mathbb{R}_{+}}\min \left[\frac{k^{2}}{16}\delta,\frac{k}{4}(
\epsilon_{1}\epsilon_{2}(1-\delta)-\delta),\right.\nonumber\\
&&\,\,\,\,\,\,\,\,\left.\left(k\left(\frac{(2-\epsilon)(1-\epsilon_1)(1-\epsilon_2)(1-\delta)}{(1-\epsilon)K}-\frac{1}{1-\epsilon}\right)-\frac{\epsilon}{1-\epsilon}\right)\epsilon\right].
\end{eqnarray*}
\end{lemma}

\proof
During the proof we use the relationship $\pi(dx)p(x,dy)=\pi(dy)p(y,dx)$ without mentioning this fact.
For $A\in\mathcal{F}$ we have
\begin{eqnarray}
k_{2}(A)&=&\frac{1}{\pi(A)\pi(A^{c})}\int_{A}\pi(dx)p^{2}(x,A^{c})\nonumber \\
&=&\frac{1}{\pi(A)\pi(A^{c})}\left(\int_{A}\pi(dx)\int_{A}p(x,dy)p(y,A^{c})+\int_{A}\pi(dx)\int_{A^{c}}p(x,dy)p(y,A^{c})\right)\nonumber\\
&=&\frac{1}{\pi(A)\pi(A^{c})}\left(\int_{A}\pi(dx)\int_{A}p(x,dy)p(y,A^{c})+\int_{A^{c}}\pi(dx)\int_{A^{c}}p(x,dy)p(y,A)\right)\nonumber\\
&\ge&\inf_{A\in\mathcal{F}:\pi(A)\le\frac{1}{2}}\frac{1}{\pi(A)\pi(A^{c})}\int_{A}\pi(dx)\int_{A}p(x,dy)p(y,A^{c})\nonumber
\end{eqnarray}
and therefore
\begin{equation}\label{nummervergeben}
k_2\ge \inf_{A\in\mathcal{F}:\pi(A)\le\frac{1}{2}}\frac{1}{\pi(A)\pi(A^{c})}\int_{A}\pi(dx)\int_{A}p(x,dy)p(y,A^{c}).
\end{equation}
Hence we can assume without loss of generality that $\pi(A)\le \frac{1}{2}$.
Let us define
\[A_{\frac{k}{4}}:=\{y\in A:p(y,A^{c})\ge\frac{k}{4}\}.\]
%It is not difficult to prove that $\pi(A_{\frac{k}{4}})\ge \frac{k}{4}\pi(A)$.
Then we obtain
\begin{equation}\label{jakob1}
\int_{A}\pi(dx)\int_{A}p(x,dy)p(y,A^{c})\ge\int_{A}\pi(dx)\int_{A_{\frac{k}{4}}}p(x,dy)p(y,A^{c})\ge\frac{k}{4}\int_{A}\pi(dx)p(x,A_{\frac{k}{4}}).
\end{equation}
Let us continue by estimating $\int_{A}\pi(dx)p(x,A_{\frac{k}{4}})$. For this reason define
\[C:=A_{\frac{k}{4}}^{c}\cap A,\]
\[\tilde{A}_{\epsilon}:=\{x\in C:p(x,A_{\frac{k}{4}})\ge\epsilon\}.\]
We now consider three cases: \\
First assume that there exists $\delta_{A}>0$, such that \begin{equation}
\pi(C)\ge\delta_{A}\pi(A).
\end{equation}
%and let $\delta_{A}$ be maximal with this property.
Then it follows:
\begin{eqnarray}
k&\le&\frac{1}{\pi(C)\pi(C^{c})}\int_{C}p(x,C^{c})\pi(dx)\nonumber\\
&=&\frac{1}{\pi(C)\pi(C^{c})}\int_{C}p(x,A_{\frac{k}{4}})\pi(dx)+\int_{C}p(x,A^{c})\pi(dx)\nonumber\\
&\le&\frac{1}{\pi(C)\pi(C^{c})}\left(\int_{\tilde{A}_{\epsilon}}p(x,A_{\frac{k}{4}})\pi(dx)+\epsilon\pi(C\cap\tilde{A}_{\epsilon}^{c}) +\frac{k}{4}\pi(C)\right)\nonumber\\
&\le& \frac{1}{\pi(C)\pi(C^{c})}\int_{\tilde{A}_{\epsilon}}p(x,A_{\frac{k}{4}})\pi(dx) +2\epsilon+\frac{k}{2}.\nonumber
\end{eqnarray}
Choosing $\epsilon=\frac{k}{8}$ we get
\[ \frac{1}{\pi(C)\pi(C^{c})}\int_{\tilde{A}_{\frac{k}{8}}}p(x,A_{\frac{k}{4}})\pi(dx)\ge\frac{k}{4}\]
and therefore
\[\int_{\tilde{A}_{\frac{k}{8}}}p(x,A_{\frac{k}{4}})\pi(dx)\ge\frac{k}{4}\pi(C)\pi(C^{c})\ge\frac{k}{4}\delta_{A}\pi(A)\pi(A^{c}).\]
Together with (\ref{jakob1}) and (\ref{nummervergeben}) this yields
\begin{equation}\label{jakob2}
k_2(A)\ge\frac{k^{2}}{16}\delta_{A}.
\end{equation}
%Since $\delta_{A}$ depends on $A$, it is possible that
%for all $A_n$ with $k_{2}(A_n)\rightarrow k_2$ the associated $\delta_{A}_{n}$ ($\delta_{A}_n$ as in the proof) converge to zero. \\

%Having this in mind, we now consider the case where
 %and $\delta_{A}>0$ minimal with this property (i.e. we allow $\pi(C)=0$).
Let us define
\[B_{\epsilon_1}:=\{x\in A_{\frac{k}{4}}:p(x,A^{c})<1-\epsilon_1\}.\]
For the second case we assume that $\pi(C)\le\delta_{A}\pi(A)$ and an $\epsilon_2> 0$ exists 
such that %maximal with the property that
\[\pi(B_{\epsilon_1})\ge\epsilon_2\pi(A_{\frac{k}{4}}).\]
With these assumptions it follows that
\begin{eqnarray}\label{jakob3}
k_2(A)&=&\frac{1}{\pi(A)\pi(A^{c})}\int_{A}\pi(dx)p^{2}(x,A^{c})\ge\frac{1}{\pi(A)\pi(A^{c})}\int_{A_{\frac{k}{4}}}\pi(dx)p^{2}(x,A^{c})\nonumber\\
&=&\frac{1}{\pi(A)\pi(A^{c})}\left(\int_{A_{\frac{k}{4}}}\pi(dx)\int_{B_{\epsilon_1}}p(x,dy)p(y,A^{c})+\int_{A_{\frac{k}{4}}}\pi(dx)\int_{B_{\epsilon_1}^{c}  }p(x,dy)p(y,A^{c})\right)\nonumber\\
&\ge&\frac{k}{4}\frac{1}{\pi(A)\pi(A^{c})}\int_{A_{\frac{k}{4}}}\pi(dx)p(x,B_{\epsilon_1})=\frac{k}{4}\frac{1}{\pi(A)\pi(A^{c})}\int_{B_{\epsilon_1}}\pi(dx)p(x,A_{\frac{k}{4}}).\nonumber
\end{eqnarray}
Moreover, we have
\begin{enumerate}
\item
\[
\int_{B_{\epsilon_1}}\pi(dx)p(x,A_{\frac{k}{4}}^{c}\cap A)=\int_{A_{\frac{k}{4}}^{c}\cap A}\pi(dx)p(x,{B_{\epsilon_1}})\le\pi(A_{\frac{k}{4}}^{c}\cap A)\le\delta_{A}\pi(A).
\]
\item
\[
\int_{B_{\epsilon_1}}\pi(dx)p(x,A)\ge\epsilon_1\pi(B_{\epsilon_1})\ge\epsilon_1\epsilon_2\pi(A_{\frac{k}{4}})\ge\epsilon_1\epsilon_2 (1-\delta_{A})\pi(A).
\]
\end{enumerate}
Subtracting the first inequality from the second we obtain
\[\int_{B_{\epsilon_1}}\pi(dx)p(x,A_{\frac{k}{4}})\ge (\epsilon_2\epsilon_1(1-\delta_{A})-\delta_{A})\pi(A)\]
and hence
\begin{equation}
k_{2}(A)\ge\frac{k}{4}(\epsilon_2\epsilon_1(1-\delta_{A})-\delta_{A}).
\end{equation}
%Diese Absch"atzung gen"ugt dann f"ur den Beweis der Aussage des lemmas, wenn eine Folge $A_n$ existiert mit $k_{2}(A_n)\rightarrow k_2$
%und wir eine von Null weg beschr"ankte Folge $\epsilon_{1}^{(n)}$ finden k"onnen, sodass die
%zugeh"orige Folge $\epsilon_{2}^{(n)}$ auch von Null weg beschr"ankt bleibt.\\
For the last case let us assume that there exists $\delta_{A}$, $\epsilon_{1}$ and $\epsilon_2$ with the same $\delta_{A}$, $\epsilon_{1}$ and $\epsilon_2$ as before such that

\[\pi(C)\le\delta_{A}\pi(A);\,\,\,\,\pi(B_{\epsilon_1})\le \epsilon_{2}\pi(A_{\frac{k}{4}}).\]
First, we observe that for small $\epsilon_1,\epsilon_2$ and $\delta_{A}$, the associated
$\pi(A)$ is bounded away from $\frac{1}{2}$. This can be established in the following way:
\begin{eqnarray}
\int_{A}\pi(dx)p(x,A^{c})&\ge&\int_{A_{\frac{k}{4}}}\pi(dx)p(x,A^{c})\ge\int_{A_{\frac{k}{4}}\cap B_{\epsilon_1}^{c}}\pi(dx)p(x,A^{c})\nonumber\\
&\ge&(1-\epsilon_1)\pi(A_{\frac{k}{4}}\cap B_{\epsilon_1}^{c})\ge(1-\epsilon_1)(1-\epsilon_2)\pi(A_{\frac{k}{4}})\nonumber\\
&\ge&(1-\epsilon_1)(1-\epsilon_2)(1-\delta_{A})\pi(A)\nonumber.
\end{eqnarray}
This yields
\[k(A)\ge\frac{(1-\epsilon_1)(1-\epsilon_2)(1-\delta_{A})}{\pi(A^{c})}.\]
Since $k(A)\le K$ for all $A\in\mathcal{F}$ by definition of $K$, we obtain
\begin{equation}\label{jakob3.5}
\pi(A^{c})\ge\frac{(1-\epsilon_1)(1-\epsilon_2)(1-\delta_{A})}{K}.
\end{equation}
This inequality will now be used in order to continue the estimation of $k_2(A)$. Define
\[H_{\epsilon}:=\{y\in A^{c}:p(y,A^{c})\ge\epsilon\}.\]
Then we have
\begin{eqnarray}\label{jakob4}
k_2(A)&\ge&\frac{1}{\pi(A)\pi(A^{c})}\int_{A}\pi(dx)\int_{A^{c}}p(x,dy)p(y,A^{c})\nonumber\\
&\ge&\frac{1}{\pi(A)\pi(A^{c})}\int_{A}\pi(dx)\int_{H_{\epsilon}}p(x,dy)p(y,A^{c})\nonumber\\
&\ge&\epsilon\underbrace{\frac{1}{\pi(A)\pi(A^{c})}\int_{A}\pi(dx)p(x,H_{\epsilon})}_{L}.
\end{eqnarray}
In order to obtain a suitable estimation of $L$, we consider the probability flow out of the set $A\cup H_{\epsilon}^{c}$:
\begin{eqnarray}
k&\le& \frac{1}{\pi(A\cup H_{\epsilon}^{c})\pi(H_{\epsilon})}\int_{A\cup H_{\epsilon}^{c}}\pi(dx)p(x,H_{\epsilon})\nonumber\\
&=&\frac{1}{\pi(A\cup H_{\epsilon}^{c})\pi(H_{\epsilon})}\left(\int_{A}\pi(dx)p(x,H_{\epsilon})+\int_{H_{\epsilon}^{c}\cap A^{c}}\pi(dx)p(x,H_{\epsilon})\right)\nonumber\\
&\le&\frac{\pi(A)\pi(A^{c})}{\pi(A\cup H_{\epsilon}^{c})\pi(H_{\epsilon})}L + \frac{\pi(H_{\epsilon}^{c}\cap A^{c})}{\pi(A\cup H_{\epsilon}^{c})\pi(H_{\epsilon})}\epsilon\nonumber\\
&\le&\frac{1}{\pi(H_{\epsilon})}\left(L+\frac{\pi(H_{\epsilon}^{c}\cap A^{c})}{\pi(A)}\epsilon\right).
\end{eqnarray}
From this we obtain
\begin{equation}\label{jakob5}
L\ge\pi(H_{\epsilon})k-\frac{\pi(H_{\epsilon}^{c}\cap A^{c})}{\pi(A)}\epsilon.
\end{equation}
Now we continue to estimate $\pi(H_{\epsilon})$ and $\frac{\pi(H_{\epsilon}^{c}\cap A^{c})}{\pi(A)}$.
It holds
\[\pi(A)\ge\int_{H_{\epsilon}^{c}\cap A^{c}}\pi(dx)p(x,A)\ge(1-\epsilon)\pi(H_{\epsilon}^{c}\cap A^{c}).\]
It follows that
\begin{equation}\label{jakob6}
\frac{\pi(H_{\epsilon}^{c}\cap A^{c})}{\pi(A)}\le\frac{1}{1-\epsilon}.
\end{equation}
By adding $\pi(H_{\epsilon})$ to both sides of $\pi(H_{\epsilon}^{c}\cap A^{c})\le\frac{1}{1-\epsilon}\pi(A)$
we obtain
\[\pi(A^{c})\le\frac{1}{1-\epsilon}\pi(A)+\pi(H_{\epsilon})\]
and hence
\[\pi(H_{\epsilon})\ge \pi(A^{c})\frac{2-\epsilon}{1-\epsilon}-\frac{1}{1-\epsilon}.\]
Using (\ref{jakob3.5}) we get
\[\pi(H_{\epsilon})\ge\frac{(1-\epsilon_1)(1-\epsilon_{2})(1-\delta_{A})}{K}\frac{2-\epsilon}{1-\epsilon} -\frac{1}{1-\epsilon}.\]
This and (\ref{jakob6}) inserted into (\ref{jakob5}) yields
\begin{equation}
L\ge k\left(\frac{(2-\epsilon)(1-\epsilon_1)(1-\epsilon_2)(1-\delta_{A})}{(1-\epsilon)K}-\frac{1}{1-\epsilon}\right)-\frac{\epsilon}{1-\epsilon}.
\end{equation}
Again, inserting this into (\ref{jakob4}), we obtain
\begin{equation}
k_{2}(A)\ge\epsilon\left(k\left(\frac{(2-\epsilon)(1-\epsilon_1)(1-\epsilon_2)(1-\delta_{A})}{(1-\epsilon)K}-\frac{1}{1-\epsilon}\right)-\frac{\epsilon}{1-\epsilon}\right).
\end{equation}
So we proved three different inequalities for $k_{2}(A)$, which are listed below:

\begin{enumerate}
\item
\[k_{2}(A)\ge\frac{k^{2}}{16}\delta_{A} \,\,\,\mbox{for }\pi(C)\ge\delta_{A}\pi(A).\]
\item
\[k_{2}(A)\ge\frac{k}{4}\left(\epsilon_1 \epsilon_2 (1-\delta_{A})-\delta_{A}\right)\,\,\,\mbox{for } \pi(C)\le \delta_{A} \pi(A),\,\,\,\,\pi(B_{\epsilon_1})\ge\epsilon_2\pi(A_{\frac{k}{4}}).\]
\item
\begin{eqnarray}
k_{2}(A)&\ge&\epsilon\left(k\left(\frac{(2-\epsilon)(1-\epsilon_1)(1-\epsilon_2)(1-\delta_{A})}{(1-\epsilon)K}-\frac{1}{1-\epsilon}\right)-\frac{\epsilon}{1-\epsilon}\right)\nonumber\\
&&\mbox{for }\pi(C)\le\delta_{A}\pi(A),\,\,\,\, \pi(B_{\epsilon_1})\le\epsilon_2\pi(A_{\frac{k}{4}}).\nonumber\\
\end{eqnarray}
\end{enumerate}
As an immediate consequence we obtain
\begin{eqnarray}\label{wua}
k_2&\ge&\sup_{\delta,\epsilon_{1},\epsilon_{2},\epsilon\in\Rset_{+}}\min \left[\frac{k^{2}}{16}\delta,\frac{k}{4}(
\epsilon_{1}\epsilon_{2}(1-\delta)-\delta),\right.\nonumber\\
&&\,\,\,\,\,\,\,\,\left.\left(k\left(\frac{(2-\epsilon)(1-\epsilon_1)(1-\epsilon_2)(1-\delta)}{(1-\epsilon)K}-\frac{1}{1-\epsilon}\right)-\frac{\epsilon}{1-\epsilon}\right)\epsilon\right].
\end{eqnarray}
This finishes the proof of the lemma.
\finishproof
\begin{lemma}\label{ja7}
Let $\xi_1,\xi_2,\ldots$ be a reversible Markov chain such that
\[0<k=\inf_{A\in\mathcal{F}}k(A)\le \sup_{A\in\mathcal{F}}k(A):=K<2.\] 
Then it follows that
\[k_2>0.\]
\end{lemma}
\proof
Choose $\epsilon_1$ and $\epsilon_2$ in a way such that $\epsilon_1\epsilon_2(1-\delta)=2\delta$  (e.g.
$\epsilon_1=\epsilon_2=\sqrt{\frac{2\delta}{1-\delta}}$). Now choose $\epsilon$ and $\delta$ sufficiently small such that the third term of (\ref{wua}) is larger than zero. 
But then (\ref{wua}) is also bounded away from zero. So the lemma is proven.
\finishproof
%Es wird nun nicht versucht, das Supremum dieses Ausdrucks zu bestimmen, da das sehr kompliziert
%erscheint (dies kann jedoch mit Hilfe eines Computers gegebenenfalls sehr gut approximiert
%werden). Es ist gelungen, einen expliziten Ausdruck in Abh"angigkeit von $k$ und $l$ anzugeben.
%Dieser ist jedoch sehr unsch"on und vermutlich weit vom Optimum entfernt und soll deshalb nicht
%angegeben werden. Vielmehr soll kurz beschrieben werden, warum obiger Ausdruck von Null weg
%beschr"ankt ist.

\begin{prop}\label{maja}
For the reversible Markov chain $\xi_1,\xi_2,\ldots$ the following holds true:
\begin{equation}\label{prelaststep}
k_2>0\,\,\Longleftrightarrow\,\,0<k\le K<2
\end{equation}
\end{prop}
\proof 
The proof follows immediately from Lemma \ref{leema} and Lemma \ref{ja7} 
\finishproof
The proof of Theorem 2 follows from Corollary \ref{majaja} and 
Proposition \ref{maja}.

\section{Acknowledgments}
The author thanks Manfred Denker for helpful conversations, Sachar Kablutschko and Jörn Dannemann for reading the preprint.

\end{document}